\documentclass[12pt,a4paper,reqno]{amsart}
\usepackage[margin=2.65cm]{geometry}
\usepackage{graphics, amsmath,amssymb, enumitem, comment, url}
\usepackage{graphicx}
\usepackage{epsfig}

\newif\ifcolorcomments
\newcommand{\allowcomments}[4]{
\newcommand{#1}[1]{\ifdraft{\ifcolorcomments{\textcolor{#4}{##1 --#3}}\else{\textsl{ ##1 \ --#3}}\fi}\else{}\fi}
}

\allowcomments{\commumtaz}{MH}{Mumtaz}{green}
\allowcomments{\comjohannes}{JS}{Johannes}{blue}
\allowcomments{\comdavid}{DS}{DS}{magenta}

\colorcommentstrue
\usepackage{amssymb}
\usepackage{amsmath,amsthm}

\usepackage{colortbl}

\usepackage{amssymb,amsmath}
\usepackage{geometry}    
\usepackage{mathrsfs}

\newtheorem{theorem}{Theorem}[section]

\theoremstyle{definition}
\newtheorem{definition}[theorem]{Definition}

\newcommand{\HD}{{\dim_H}}

\newcommand{\HH}{\mathcal H}

\newcommand{\N}{\mathbb N}% by default, $\N$ is American naturals $\{1,2,\ldots\}$
\newcommand{\Q}{\mathbb Q}

\newcommand{\R}{\mathbb R}

\newcommand{\Z}{\mathbb Z}

\newcommand{\mbf}{\mathbf}
\newcommand{\pp}{\mbf p}
\newcommand{\qq}{\mbf q}
\newcommand{\rr}{\mbf r}

\newcommand{\xx}{\mbf x}
\newcommand{\yy}{\mbf y}

\renewcommand{\text}{\textup}

\newcommand{\NPC}[1]{\ignorespaces}

\newif\ifdraft\drafttrue

\IfFileExists{marvosym.sty}{
\RequirePackage{marvosym} 
}{

}

\IfFileExists{wasysym.sty}{
\RequirePackage{wasysym}
}{

}

\let\6\sectionsymbol

\begin{document}
\title{Metrical theorems for unconventional height functions}
\author[M. Hussain]{Mumtaz ~Hussain}
\address{Mumtaz Hussain,  Department of Mathematical and Physical Sciences, La Trobe University, Bendigo 3552, Australia. }\email{m.hussain@latrobe.edu.au}
%\authordavid

%\subjclass[2010]{Primary }
%\keywords{}
%\date{}
%\dedicatory{}

\begin{abstract}  In this paper we consider the simultaneous approximation of real points by rational points with the error of approximation given by the functions of `non-standard' heights. We prove analogues of Khintchine and Jarn\'ik-Besicovitch theorems for this setting,  thus answering some questions raised by Fishman and Simmons (2017).   

\end{abstract}
\maketitle

\section{Introduction and Results}
Throughout this paper, let $d\geq 1$ be a fixed integer, $\pp=(p_1,\ldots, p_d)\in \Z^d$,  and $\xx=(x_1,\dots,x_d)\in\R^d$. 
Let $\psi:\N\to[0, \infty)$ be an \emph{approximating function} such that $\psi(n)\to 0$  as $n\to \infty$.
A central theme in simultaneous Diophantine approximation is to estimate the `size'  of the set
$$S_d(\psi):=\{\xx=(x_1, \ldots, x_d)\in \R^d: \|\xx-\pp/q\|<\psi(q) \ \text{for infinitely many} \ \pp/q\in\Q^d\},$$
where $\|\cdot\|$ denote the max norm.  Here, the approximating function $\psi$ is a function of  the `standard height' of the rational point $\pp/q$, that is,  $H(\pp/q)=q$.    The metrical theory (Dirichlet, Khintchine, Jarn\'ik-Besicovitch, and Jarn\'ik theorems) for this set is well-known,  see \cite{BeresnevichVelani} for further details.  However, if the height function $H:\Q^d\to\N$ is nonstandard then the metrical theory is not so well developed. In this direction, the best known result is by Fishman-Simmons \cite{FishmanSimmons3} calculating the irrationality exponents for every point $\R^d\setminus\Q^d$ approximated by rationals with certain nonstandard height functions. 

As in \cite{FishmanSimmons3}, let $\Theta:\N^d\to \N$ and let $H_\Theta:\Q^d\to\N$ be defined by the formula $$H_\Theta\left(\frac {p_1}{q_1}, \ldots, \frac{p_d}{q_d} \right)=\Theta(q_1, \ldots, q_d).$$ Here the rationals $p_i/q_i$ are in reduced form. The set $S_d(\psi)$ is concerned with the standard height function $H_{\rm lcm}$, where ${\rm lcm} : \N^d \to \N$ is the least common multiple function. That is, given $\rr\in\Q^d$, $H_{\rm lcm}(\rr)$ is the smallest number $q$ such that $\rr=\frac\pp q$ for some $\pp\in\Z^d$.  As in \cite{FishmanSimmons3}, various heights can be given as \begin{align*}
H_{\max}(\rr) &=\max (q_1, \ldots, q_d), \\
 H_{\min}(\rr) &=\min (q_1, \ldots, q_d), \\
 H_{\rm prod}(\rr) &=q_1\times \ldots\times q_d.
 \end{align*} 
For convenience we will denote $\psi_\tau(q):=q^{-\tau}$ for some $\tau\geq 0$.  Note that the height functions mentioned above are only different if $d>1$. For $d=1$ the height of the rational $p/q$ will be $q$ and therefore  $H_{\max}(p/q) = H_{\min}(p/q) = H_{\rm prod}(p/q) = H_{\rm lcm}(p/q) = q$.
\begin{definition}Given a height function $H : \Q^d \to \N$ and a point $\xx\in\R^d\setminus\Q^d$, the {\em exponent of irrationality} of $\xx$ is defined as
\begin{equation*}\label{exp}
\omega_H(\xx)=\liminf_{\substack{ {\rr\in\Q^d}\\{\rr\to \xx} }} \frac{-\log\|\xx-\rr\|}{\log H(\rr)}=\lim_{\epsilon\to 0} \inf_{\substack{ {\rr\in\Q^d}\\ {\|\xx-\rr\|\leq \epsilon }}}\frac{-\log\|\xx-\rr\|}{\log H(\rr)}
\end{equation*} 
Equivalently\footnote{To keep notation simple and keeping length of this note short, we will not go in detailed description of notations and refer the reader to Fishman-Simmons \cite{FishmanSimmons3} for details.}, $\omega_H(\xx)$  is the supremum of  $\tau\geq 0$ such that
$$\|\xx-\rr_n\|<\psi_\tau\circ H(\rr_n)\quad \text{for some sequence} \ \Q^d\ni \rr_n\to \xx$$
The exponent of irrationality of the height function $H$ is the number
$$\omega_d(H)=\inf_{\xx\in \R^d\setminus Q^d}\omega_H(\xx).$$
\end{definition}
Notice that Dirichlet's theorem implies that $\omega_d(H_{\rm lcm})= 1+1/d$.

\begin{definition}Given a height function $H : \Q^d \to \N$, a function $\psi:\N \to (0, \infty)$, and a point $\xx\in\R^d$, let
\begin{equation*}\label{chpsi}
C_{H, \psi}(\xx)=\liminf_{\substack{ {\rr\in\Q^d}\\{\rr\to \xx} }}\frac{\|\xx-\rr\|}{\psi\circ H(\rr)}
\end{equation*}
 
Equivalently, $C_{H, \psi}(\xx)$  is the infimum of all $C\geq 0$ such that
$$\|\xx-\rr_n\|<C\psi\circ H(\rr)\quad \text{for some sequence} \ \Q^d\ni \rr_n\to \xx$$

\end{definition}
The main result of \cite[Theorem 1.1]{FishmanSimmons3} is the following theorem.

\begin{theorem}[Fishman-Simmons, \cite{FishmanSimmons3}]\label{FStheorem} The exponents of irrationality of $H_{\max}$, $H_{\min}$ and $H_{\rm prod}$ are

\begin{align*}
\omega_d\left(H_{\max}\right)&=\frac{d}{(d-1)^{(d-1)/d} } \quad if \ d\geq 2\\
\omega_d\left(H_{\min}\right) &=2\\
\omega_d\left(H_{\rm prod}\right) &=\frac{2}{d}.
\end{align*}

\end{theorem}
In other words, Theorem \ref{FStheorem} concerns exponents valid `everywhere',  that is,  functions $\psi$ for which $C_{H, \psi} (\xx)<\infty$ for every point 
$\xx\in\R^d\setminus\Q^d$. At the end of their paper, Fishman--Simmons posed the question, 

{\em ``what happens  if ``every''  is replaced by ``almost every''-- with respect to Lebesgue measure or even with respect to some fractal measure? 
 Once we know what
``almost every''  point does, it can be asked what is the Hausdorff dimension of the set
of exceptions, i.e. the set of $\xx$ which behave differently from almost every point.''} 

We answer these questions below.

\begin{theorem}
\label{theoremkhinchin}
Almost any point $\xx\in\R^d$ satisfies
\[
\omega_{H_{\max}}(\xx) = \omega_{H_{\rm prod^{1/d}}}(\xx) = \omega_{H_{\min}}(\xx) = 2.
\]
\end{theorem}
Comparing Theorems \ref{FStheorem} and \ref{theoremkhinchin}, we note that the exponents for $H_{\max}$ and $H_{\rm prod}$ are different for $d>2$. It raises a natural question of size of the exceptions in terms of Hausdorff dimension. For a set $X$, let $\HD(X)$ denote the Hausdorff dimension of the set $X$. For the definition of Hausdorff dimension and measure we refer the reader to \cite{BernikDodson}.

Next, define the set
\begin{equation*}\label{WPhi}
S_\Theta(\tau) = \left\{\xx\in \R^d : \begin{array}{ll} \text{ the system $|x_i - p_i/q_i| < \Theta(\qq)^{-\tau} \;\; (i=1,\ldots,d)$}\\ [1ex] \text{ has infinitely many solutions $(p_1/q_1,\ldots,p_d/q_d)\in\Q^d$} \end{array}\right\}.
\end{equation*}

\begin{theorem}\label{thmJB}
Let $\Theta \in \{\max, \  {\rm prod}^{1/d}\}$. For all $\tau\geq 2$,
\[
\HD(S_\Theta(\tau)) = \frac{2d}{\tau}\cdot
\]
\end{theorem}

In the case that $\Theta=\min$, we prove the following result by first tailoring the set $S_d( \psi)$ according to the $\min$ height function. Define,
\begin{equation*}\label{Wmin}
S_d(\psi, \min) =\left\{\xx\in \R^d : \begin{array}{ll} \text{ the system $|x_i - p_i/q_i| < \psi({\rm min}(\qq)) \;\; (i=1,\ldots,d)$}\\ [1ex] \text{ has infinitely many solutions $(p_1/q_1,\ldots,p_d/q_d)\in\Q^d$} \end{array}\right\}.
\end{equation*}

\begin{theorem}\label{Smintheorem}
For all monotonic approximating functions $\psi$, if $\Theta=\min$, then we have
\[
S_d(\psi,\min) = \bigcup_{i=1}^d \left(\R^{i-1} \times S_1(\psi)\times \R^{d-i}\right).
\]
\end{theorem}

\section{Proofs}

\noindent{\bf Notation.} By $B:=B(\xx, \rho)=\{\yy\in\R^d: |\yy-\xx|\leq \rho\}$ we mean the ball centred at the point $\xx\in\mathbb{R}^d$ of radius $\rho$. Define another ball $B^s:=B(\xx, \rho^{s/d})$.
For real quantities $A,B$ and a parameter $t$, we write $A \lesssim_t B$ if $A \leq c(t) B$ for a constant $c(t) > 0$ that depends on $t$ only (while $A$ and $B$ may depend on other parameters). 
We write  $A\asymp_{t} B$ if $A\lesssim_{t} B\lesssim_{t} A$.
If the constant $c>0$ depends only on parameters that are constant throughout a proof, we simply write $A\lesssim B$ and $B\asymp A$.

\subsection{Mass Transference Principle}

In a landmark paper \cite{BeresnevichVelani}, Beresnevich-Velani  introduced the Mass Transference Principle which has become a major tool in transferring the Lebesgue measure theoretic statements for limsup sets defined by balls to the Hausdorff measure statements for limsup sets defined by balls. This is surprising as the Lebesgue measure is the `coarser' notion of `size'  than the Hausdorff measure.  Although this principle holds for general metric spaces with distance induced by any fixed norm (not necessarily Euclidean), we  only need the following version suitable for the purposes of this paper. 
%By a ball $B:=B(\xx, \rho)\subseteq \R^d$ centred at $\xx\in\R^d$ and radius $\rho$ we mean  the set of all $\yy\in\R^d$ such that $|\xx-\yy|\leq r$. 
Let $\{B(x_n, \rho_n)\}_{n\in\N}$ be a sequence of balls in $\R^d$ with $\rho_n\to 0$ as $n\to \infty$. 

\begin{theorem}[Beresnevich-Velani, 2006]\label{mtp} Let $s>0$ and suppose for any  ball $B$ in $\R^d$, 
\begin{equation*}\label{eq1mtp}\HH^d(B\cap\limsup_{i\to\infty}B_i^s)=\HH^d(B).\end{equation*}
 Then for any ball $B\subset X$
\begin{equation*}\label{eq2mtp}\HH^s(B\cap\limsup_{i\to\infty}B^d_i)=\HH^s(B).\end{equation*}
\end{theorem}
Here note that $\HH^d$ is simply the $d$-dimensional Lebesgue measure. Roughly speaking, given that a limsup set defined by balls has full Lebesgue measure, the Hausdorff measure of the limsup set with the radius of each generating balls shrinking in all directions at the same rate can be calculated.  As one would expect the Mass Transference Principle has many applications in number theory such as the derivation of Jarn\'ik's theorem from Khintchine's theorem or the implication of Dirichlet's theorem to the Jarn\'ik-Besicovitich theorem, see sections 3.1 and 3.2 in \cite{BeresnevichVelani}. There have been several generalisation of this fundamental result, see for instance \cite{ AllenBaker, HussainSimmons4, WangWu19}.

\subsection{Continued Fractions}
Every irrational $x\in [ 0,1)$ can be uniquely expressed as a simple
infinite continued fraction expansion of the form: 
\begin{equation*}
x=[a_{1}(x),a_{2}(x),\ldots ,],
\end{equation*}
where $a_{n}(x)\in \mathbb{N},$ $n\geq 1$, are known as the \emph{partial quotients}
of $x.$ We denote $q_n=q_n(x)$ to be the denominator of the $n$th convergent to $x$, that is,
$$[a_{1}(x),a_{2}(x),\ldots, a_n(x)]:=\frac{p_n(x)}{q_n(x)}.$$

The following recursive relations are well-known
\begin{equation*}
\begin{split}
p_{-1}& =1,~p_{0}=0,~p_{n+1}=a_{n+1}p_{n}+p_{n-1}, \\
q_{-1}& =0,~q_{0}=1,~q_{n+1}=a_{n+1}q_{n}+q_{n-1},
\end{split}
\label{recu}
\end{equation*}
 and then, for any $n\geq 1$, we have
\begin{equation}
\label{qnqn1}
\frac{1}{3a_{n+1}q_{n}^{2}}\,<\,\Big|x-\frac{p_{n}}{q_{n}}\Big| \asymp \frac{1}{q_{n+1}q_{n}}.
\end{equation}

Fix any $\xx\in\R^d\setminus\Q^d$, and for each $i=1, \ldots, d,$ let $\left(\frac{p_n^{(i)}}{q_n^{(i)}} \right)_{n=1}^{N_i}$ be the convergents of the continued fraction expansion of $x_i$. Here $N_i\in\N\cup\{\infty\}$ with $N_i=\infty$ for at least one $i$.

\subsection{Proof of Theorem \ref{theoremkhinchin}} Denote $q_{k,i} = q_k(x_i)$ for any $ i= 1,\ldots,d$. Then, it is straightforward to see that almost any point $\xx$ satisfies
\[
q_{k+1,i} \lesssim q_{k,i}^{1+\epsilon}
\]
for all $\epsilon > 0$.  Fix $\Theta\in\{\max, \min, {\rm prod}^{1/d}\}$ and if $\omega > \omega_{H_\Theta}(\xx)$, then
\[
\|\xx - \rr\| =\max_{i=1}^d|x_i-r_i|\gtrsim H_\Theta(\rr)^{-\omega} = \Theta(\qq)^{-\omega}.
\]
Thus by \eqref{qnqn1},
\begin{align*}
\max_i \frac{1}{q_{n_i,i} q_{n_i+1,i}} &\gtrsim \Theta(\qq)^{-\omega}\\
\min_i [\log q_{n_i,i} + \log q_{n_i+1,i}] &\lesssim_+ \omega \log \Theta(\qq)\\
\min_i (2+o(1))\log q_{n_i,i} &\lesssim_+ \omega\Theta((\log q_{n_i,i})_{i=1}^d) \\ &\leq \omega \max_i(\log q_{n_i,i}).
\end{align*}
Here $\Theta(\log q_{n_i,i})$ is the point whose $i$th coordinate is $\log(q_{n_i,i})$. Next, by choosing $n_i$ such that $$\log q_{n_i,i} = (1 + o(1)) \log q_{n_1,1},$$ we get $\omega \geq 2$. The converse is similar. 
%{\color{red}(Simon - What is the $q_{n_j,j}$ here? Where did the $j$ subscripts come from? )} \comdavid{$j$ can be anything, I changed it to 1}
%\end{proof}

\subsection{Proof of Theorem \ref{thmJB}} 
For the lower bound, we use the mass transference principle Theorem \ref{mtp}.  To this end, we fix $\epsilon>0$ and note that by Theorem \ref{theoremkhinchin}, for almost any point $\xx\in \R^d$, there are infinitely many rational solutions $(p_1/q_1,\ldots,p_d/q_d)\in\Q^d$ to the inequality
$$\|x_i - p_i/q_i\| < \Theta(\qq)^{-(2-\epsilon)}$$
This implies that $$\HH^d(S_\Theta(2-\epsilon)\cap [0, 1)^d)=1.$$ Thus, by Theorem \ref{mtp},  with $s={\frac{(2-\epsilon)d}{\tau}}$,
%\quad \phi_i(q) = q^{-(2-\epsilon)},\quad \psi_i(q) = q^{-\tau},$
we have that
$$\HH^{\frac{(2-\epsilon)d}{\tau}}(S_\Theta(\tau))=\infty.$$
Since $\epsilon$ was arbitrary we get $\HD(S_\Theta(\tau)) \geq \frac{2d}{\tau}$. 

For the upper bound, we construct a simple cover of the set $S_\Theta(\tau)$ and then use the Hausdorff-Cantelli Lemma \cite{BernikDodson}. Recall that, for the set
\[S_\Theta(\tau)=\limsup_{\Theta(\qq)\to\infty}\left\{\xx\in\R^d: \|x_i - p_i/q_i\| < \Theta(\qq)^{-\tau}\right\},\]
there are two cases: 

\noindent{\bf Case I.} If $\Theta(\qq)=\max(q_1, \ldots, q_d)=q$, then the $s$-dimensional Hausdorff measure of $S_\Theta(\tau)$ becomes
\begin{equation}\label{eqcaseI}
\HH^s\left(S_\Theta(\tau)\right)\ll \sum_{q}q^{2d-1-\tau s}.
\end{equation}
\noindent{\bf Case II.} If $\Theta(\qq)={\rm prod}^{\frac1d}=(q_1\times \ldots\times q_d)^{\frac1d}$, then
\begin{align}\label{eqcaseII}
\HH^s\left(S_\Theta(\tau)\right)&\ll \sum_{q_1,\ldots,q_d} q_1 \cdots q_d (q_1\ldots q_d)^{-\tau s/d} = \left(\sum_q q^{1-\tau s/d}\right)^d.
\end{align}
It is straightforward to see that the series appearing in the right hand side of both the equations \eqref{eqcaseI} and \eqref{eqcaseII} converges if $s > \frac{2d}{\tau}$. Using the Hausdorff-Cantelli Lemma \cite{BernikDodson}, we have that  $\HH^s\left(S_\Theta(\tau)\right)=0$ for any $s > \frac{2d}{\tau}$. Hence $\HD(S_\Theta(\tau))  \leq \frac{2d}{\tau}.$

\subsection{Proof of Theorem \ref{Smintheorem}} 
The proof is rather straightforward. Let $\xx\in S_d(\psi,\min)$, then there exists infinitely many $\rr$ such that $ \|\xx - \rr\| \leq \psi(\min(\qq))$. Since $\psi$ is decreasing, we have $$\psi(\min(\qq))\leq \max_i(\psi(q_i)).$$ Therefore, 
\begin{align*}
\xx\in S_d(\psi,\min)
%&\Longleftrightarrow\;\;  \text{there exists infinitely many} \   \rr \ \text {such that} \  \|\xx - \rr\| \leq \psi(\min(\qq))\\
%&\Longleftrightarrow\;\;  \text{there exists infinitely many} \   \rr \ \text {such that}  \|\xx - \rr\| \leq \max_i(\psi(q_i)) \quad(\dag)\\
&\Longleftrightarrow\;\;  \exists \  \text{ infinitely many} \  \rr ,\ \text {and} \  1\leq i\leq d \ \text{such that}\; \|\xx - \rr\| \leq \psi(q_i)\\
&\Longleftrightarrow\;\;  \exists\  1\leq i\leq d, \text{ and infinitely many} \   \rr \  \ \text{such that}\; \|\xx - \rr\| \leq \psi(q_i)\\
&\Longleftrightarrow\;\;  \exists\   1\leq i\leq d, \text{ and infinitely many} \   r_i \  \ \text{such that}\;  \; \|x_i - r_i\| \leq \psi(q_i)\quad(\dag)\\
&\Longleftrightarrow\;\;  \text {there exists}  \  1\leq i\leq d, \ \text{such that}  \; x_i \in S_1(\psi).
\end{align*}
where in the equation $(\dag)$ we use the fact that reals can be approximated arbitrarily well by rationals.
%\end{proof}

\medskip

\medskip

\noindent{\bf Acknowledgements.} This research is supported by the ARC DP200100994. The author thanks David Simmons and Simon Baker for useful discussions.

%
%%
%\bibliographystyle{amsplain}
%%
%\bibliography{../../../bibliography4}

\providecommand{\bysame}{\leavevmode\hbox to3em{\hrulefill}\thinspace}
\providecommand{\MR}{\relax\ifhmode\unskip\space\fi MR }
\providecommand{\MRhref}[2]{%
  \href{http://www.ams.org/mathscinet-getitem?mr=#1}{#2}
}
\providecommand{\href}[2]{#2}

\end{document}